\input amstex.tex
\documentstyle{amsppt}
%
%
%
%
%
\magnification\magstephalf

\def\slmn{{\frak{sl}}(M+1|N+1)}
\def\Uqslmn{U_q({\frak{sl}}(M+1|N+1))}

\def\UqMN{U_{q}({\widehat{\frak{sl}}(M+1|N+1))}}

\def\ignore#1{ } 
\NoBlackBoxes
\TagsOnRight
\topmatter
\rightline{OIT-96-1}
\rightline{Dec. 1996}
\vskip 10pt
\title
$q$-differential operator representation of
the quantum superalgebra$\Uqslmn$
\endtitle
\author
Kazuhiro Kimura
\endauthor
\affil
Department of Physics, Osaka Institute of Technology, Omiya,
Osaka 535, Japan \\
\endaffil
\abstract
A representation of the quantum superalgebra $\Uqslmn$ is constructed
based on the $q$-differential operators acting on the coherent states 
parameterized by coordinates. These coordinates correspond to
the local ones of the flag manifold. This realization provides us 
with a guide to construct the free field realization for the
quantum affine superalgebra $\UqMN$ at arbitrary level.
\endabstract
\endtopmatter
\document
\subhead 1. Introduction \endsubhead 
\par
\noindent
It has been established that affine quantum algebras play an
important role in studies on the mathematical physics and solvable
systems. It has also been cleared 
that the representation theory of these algebras,particularly,
the free field realization[1-3] provides us a powerful tools.
In order to construct the free field realization
for the quantum algebra $U_q(\hat{\frak g})$[4-9],it gives us a guide
 that the corresponding Lie algebras
$\hat {\frak g}$ is realized in terms of the $q$-differential operators
,in other words, the $q$-deformed harmonic oscillators[10-13].
 
On the other hand,  affine superalgebras have
also attracted attention in relation to superconformal algebras
and superstrings and topological field theories.
There is possibility for the $q$-deformed models of these models
to have physical meanings.
So it is worthwhile 
to investigate the representation of affine quantum superalgebras[14-17].
In the previous work[18] we have constructed a level-one representation of
the quantum affine superalgebra $\UqMN$ and the vertex operators associated 
with representations on the line of the Frenkel and Jing's construction[19].
In order to extend the representation to the case of arbitrary level, we have to 
consult the forms of the Wakimoto construction with ghost systems[1].
It is well known that there is the close relation between the Wakimoto 
construction and the realization by the differential operators on the 
flag manifold.

The aim of this paper is to construct the realization based on the 
$q$-differential operators acting on the space of analytic functions
on the coordinates of the coherent states
corresponding to the local coordinates of the flag manifold.
There have been already the papers to construct the representations
of quantum Lie superalgebras in terms of $q$-oscillators or 
$q$-differential operators[20,21]. We construct the representation
of the quantum Lie superalgebra $\Uqslmn$ on the same method in the paper[21].
This realization will
help us to construct the free field realization for the quantum affine
algebra $\UqMN$ at arbitrary level, which is the $q$-deformation
of the Wakimoto construction.
\subhead 2.
Differential operator realization of $\slmn$
\endsubhead
\par
\noindent
The Cartan matrix of the Lie superalgebra $\slmn$ is written as
$a_{ij}=(\nu_i+\nu_{i+1})\delta_{i,j}-\nu_i\delta_{i,j+1}-\nu_{i+1}\delta_{i+1,j}$
$(i,j=1,2\cdot\cdot\cdot ,M+N+1)$
,where $\nu_j=1$ for $j=1,\cdots,M+1$ and
$\nu_j=-1$ for $j=M+2,\cdots,M+N+2$.
The Lie superalgebra  $\slmn$ is defined by the Chevalley generators
$h_i,e_i,f_i(i=1,\cdots,M+N+1)$. The $\Bbb{Z}_2-\text{grading}|\cdot| \to
\Bbb{Z}_2$ of the generators are:$|e_{M+1}|=|f_{M+1}|=1$ and zero otherwise.
The relations among these generators are
$$
\align
[h_i,h_j]=&0, \tag 1  \\
[h_i,e_j]=&a_{ij}e_j, \tag 2  \\
[h_i,f_j]=&-a_{ij}f_j, \tag 3 \\
[e_i,f_j]=&\delta_{i,j}h_{ij}, \tag 4 \\
\endalign
$$
$$
\left.
\aligned
[e_j,[e_j,e_i]]&=0\\
[f_j,[f_j,f_i]]&=0
\endaligned\;\;\right\}
\qquad\text{for}\;|a_{ij}| = 1,j\not=M+1, \tag 5
$$
$$
\align
[e_{M+1},[e_{M+2},[e_{M+1},e_M]]]=&0, \tag 6 \\
[f_{M+1},[f_{M+2},[f_{M+1},f_M]]]=&0, \tag 7 \\
\endalign
$$
where we have used the notations
$[X,Y]=XY-(-1)^{|X||Y|} YX$.
For the odd generator $e_{M+1}$, the Serre relations are changed.

Let us introduce an orthonormal basis $\{\varepsilon_i|i=1,\cdots,M+N+1\}$
with the bilinear form $(\varepsilon_i|\varepsilon_j)=\nu_i\delta_{i,j}$.
The simple roots are written as $\alpha_i=\nu_i\varepsilon_i-\nu_{i+1}
\varepsilon_{i+1}$ and $(\alpha_i|\alpha_j)=a_{ij}$ is satisfied.
Let $V_\Lambda$ be the Verma module over $\slmn$ generated by the highest
weight vector $|\Lambda>$ which satisfies the highest weight conditions:
$e_i|\Lambda\!>=0,h_i|\Lambda\!>=(\alpha_i|\Lambda)|\Lambda\!>$.
Weight vectors are written as 
$$
|\Lambda;n_{lm}\!>=(X^-_{lm})^{n_{lm}}|\Lambda\!>
(1\leq l\leq m\leq M+N+1), \tag 8
$$
where 
are the generators correspond to the roots
$\alpha_{lm}=\alpha_l+\alpha_{l+1}+\dots +\alpha_m(l<m)$.
The dual module $V_\lambda^*$ is generated by $<\!\Lambda|$ such that
$<\!\Lambda|f_i=0,<\!\Lambda|h_i=(\alpha_i|\Lambda)<\!\Lambda|$,
whose dual weight vectors are
$
<\Lambda;n_{lm}|=<\Lambda;n_{lm}|(X^+_{lm})^{n_{lm}}
(1\leq l\leq m\leq M+N+1). 
$
Let's introduce the coherent state:
$$
|{\bold x}>=|x_{11},x_{12},\dots,x_{M+N+1M+N+1}>=
\prod_{l\leq m}\sum_{n_{lm}=0}^{\infty}
{x^{n_{lm}} \over n_{lm}!}
|\Lambda;n_{lm}>=\prod_{l\leq m}e^{x_{lm}X_{lm}^-}|\Lambda>. \tag 9
$$
Here $x_{lm}(1\leq l\leq m\leq M+N+1)$ are the coordinates defined as
$$
x_{lm}=
\left\{
\aligned &z_{lm}(1\leq l \leq m \leq M;M+2\leq l \leq m \leq M+N+1)\\
&\theta_{lm}(1\leq l \leq M+1,M+1\leq m\leq M+N+1), 
\endaligned 
\right. \tag 10
$$
where $z_{lm}$ are complex variables and $\theta_{lm}$ are Grassmann odd ones
which satisfy the relations $\theta_{lm} \theta_{lm}=0,\theta_{lm} \theta_{ij}=
-\theta_{ij}\theta_{lm} (l,m)\neq(i,j)$.The bilinear form $V_\Lambda^*
\otimes V_\Lambda \to {\Bbb C}$is defined by $<\!\Lambda|\Lambda\!>=1$
and $(<\!{\bold x}|X)|{\bold y}\!>=<\!{\bold x}|(X|{\bold y}\!>)$
for any $<\!{\bold x}|\in V^*_\Lambda,|{\bold y}\!>\in V_\Lambda$
and $X\in \slmn$. We can construct the representation in the space
of analytic functions on the coordinates $\{x_{ij}\}$ through 
the bilinear form.
The coordinates
\footnote{$\theta_{ij}^{-1}$ is defined as a differential operator,
$\partial_{\theta_{ij}}$} $x_{ij},x_{ij}^{-1}$
 and the differential operators
$M_{lm}=x_{lm}{\partial \over \partial x_{lm}}
(1\leq l\leq m\leq M+N+1)$ generate the Heisenberg algebra ${\Cal H}_L$
with relations:$[M_{ij},x_{lm}^\pm]=\pm\delta_{i,l}\delta_{j,m}x_{lm}^\pm$.
\proclaim{Proposition 1}
The superalgebra $\slmn$ is realized by the Heisenberg algebra 
${\Cal H}_L$.
$$
\align
h_i&=-\sum^{M+N+1}_{m=1}\sum^{m}_{l=1}\sum^{m}_{r=l}a_{ir}M_{lm}
+(\alpha_i|\Lambda) \tag 11 \\
e_i&=x_{ii}^{-1}M_{ii}+\sum_{l=1}^{i-1}x_{li-1}x_{li}^{-1}M_{li} \tag 12 \\
f_i&=\nu_{i}\sum_{l=1}^{i-1}x_{li}x_{li-1}^{-1}M_{li-1}
-\nu_{i+1}\sum_{l=i+1}^{M+N+1}x_{il}x_{i+1l}^{-1}M_{i+1l} \\
&+x_{ii}\biggl(-{a_{ii} \over 2}M_{ii}-
\sum^{M+N+1}_{m=i+1}\sum^{m}_{l=1}\sum^{m}_{r=l}a_{ir}M_{lm}
+(\alpha_i|\Lambda) \biggr) \tag 13 \\
\endalign
$$
\endproclaim
\leftline{proof.}
\par
\noindent
We can determine the form of the generators 
using the following commutation relations:
$$
\align
&[h_i,e^{x_{lm}X_{lm}^-}]=-\sum_{r=l}^{m}a_{ir}M_{rm}e^{x_{lm}X_{lm}^-}, \tag 14 \\
&[e_i,e^{x_{jj}X_{jj}^-}]|\Lambda\!>=\delta_{i,j}
x_{ii}\biggr(-{a_{ii} \over 2}M_{ii}+(\alpha_i|\Lambda)\biggr)e^{x_{ii}X_{ii}^-}
|\Lambda\!>, \tag 15 \\
&[e_i,X_{ji}]=\nu_{i}X_{ji-1}^-, \tag 16 \\
&[e_i,X_{ij}]=-\nu_{i+1}X_{i+1j}^-, \tag 17 \\
&[f_i,X_{ji-1}]=X_{ij}^-, \tag 18 \\
&[f_i,X_{ji}^-]=0,\quad [f_j,X_{ji}^-]=0,
\quad  [f_k,X_{ji}^-]=0  \quad (j+1\le k \le i-1). \tag 19 \\
\endalign
$$
We also check the results by the explicit calculation for 
the commutations of the generators.
\subhead 3.
$q$-Differential operator realization of $\Uqslmn$
\endsubhead
\par
\noindent
Now we are the position to extend the representation on the 
Heisenberg algebra to the case of the quantum Lie superalgebra $\Uqslmn$.
 If there is no confusion, we adopt the same notations in the case of $\slmn$.
The relations among the Chevalley generators
$t_i,e_i,f_i(i=1,\cdots,M+N+1)$ are written by
$$
\align
&t_it_j=t_jt_i, \tag 20 \\
&t_ie_jt_i^{-1}=q^{a_{ij}}e_j,\tag 21\\
&t_if_jt_i^{-1}=q^{-a_{ij}}e_j,\tag 22\\
&[e_i,f_j]=\delta_{i,j}{t_i-t_i^{-1} \over q-q^{-1}},\tag 23 \\
\endalign
$$
$$
\left.
\aligned
[e_j,[e_j,e_i]_{q^{-1}}]_q&=0\\
[f_j,[f_j,f_i]_{q^{-1}}]_q&=0
\endaligned\;\;\right\}
\qquad\text{for}\;|a_{ij}| = 1,j\not=M+1, \tag 24
$$
$$
\align
[e_{M+1},[e_{M+2},[e_{M+1},e_M]_{q^{-1}}]_q]=&0, \tag 25 \\
[f_{M+1},[f_{M+2},[f_{M+1},f_M]_{q^{-1}}]_q]=&0, \tag 26 \\
\endalign
$$
where we have used the notations
$[X,Y]_\xi=XY-(-1)^{|X||Y|}\xi YX$
and $q$ is complex number such that  $|q|\ne 1$.
The deformation of the differential operator is defined by
the difference operator\footnote{ There is no change in the case of Grassmann coordinates,
$D_{\theta_{ij}}=\theta_{ij}^{-1}[M_{ij}]
=\theta_{ij}^{-1}M_{ij}=\partial_{\theta{ij}}$}
 on
the coordinates$(x_{ij})
(1\leq i \leq j \leq M+N+1)$ as
$$
D_{x_{ij}}f(x)={f(qx_{ij})-f(q^{-1}x_{ij}) \over q-q^{-1}}
={q^{M_{ij}}-q^{-M_{ij}} \over x_{ij}(q-q^{-1})}f(x)
=x_{ij}^{-1}[M_{ij}]f(x),  \tag 27
$$
where $[n]$ is a standard notation ${q^n-q^{-n} \over q-q^{-1}}$.
The quantum affine super algebra $\Uqslmn$
can be endowed with the graded Hopf algebra structure.
We take the following coproduct
$$
\align
\Delta(e_i)=e_i\otimes 1+t_i\otimes e_i, \;\;
\Delta(f_i)=f_i\otimes t_i^{-1}+1\otimes f_i ,\;\;
\Delta(t_i^{\pm 1})=t_i^{\pm 1}\otimes t_i^{\pm 1},  \tag 28
\endalign
$$
and the antipode
$$
a(e_i)=-t_i^{-1} e_i,\quad a(f_i)=-f_i t_i,\quad
a(t_i^{\pm 1})=t_i^{\mp 1}. \tag 29
$$
The coproduct is an algebra automorphism
$\Delta(xy)=\Delta(x)\Delta(y)$ and the antipode is a graded
algebra anti-automorphism
$a(xy)=(-1)^{|x||y|}a(y)a(x)$ for $x,y \in \Uqslmn$.

Let  $V_\Lambda$ be the Verma module over $\Uqslmn$ generated by the 
highest weight vector $|\Lambda\!>$ such that $e_i|\Lambda\!>=0,
t_i|\Lambda\!>=q^{(\alpha_i|\Lambda)}|\Lambda\!>$.Weight vectors
are written as $|\Lambda ;n_{lm}\!>=
(X^-_{lm})^{n_{lm}}|\Lambda ;n_{lm}\!>(1\le l \le m \le M+N+1)$,
where $X^-_{ll}=f_l^-$ and $X_{lm}^-(l\le m)$ are the generators corresponding to
the non-simple roots $\alpha_l+\alpha_{l+1}+\dots+\alpha_m$:
$$
X_{lm}^-=[f_m,[f_{m-1},[f_{m-2},\dots,[f_{l+1},f_l]_{q^{-\nu_{l+1}}}
]_{q^{-\nu_{l+2}}}\dots]_{q^{-\nu_{m}}}. \tag 30
$$
The dual module $V_\Lambda^*$ is generated by $<\!\Lambda|$
such that $<\!\Lambda|f_i=0,<\!\Lambda|t_i=q^{(\alpha_i|\Lambda)}<\!\Lambda|$.
A $q$-deformed coherent state is defined by
$$
|{\bold x}>_q=|x_{11},x_{12},\dots,x_{M+N+1M+N+1}>_q=
\prod_{l\leq m}\sum_{n_{lm}=0}^{\infty}
{x^{n_{lm}} \over [n_{lm}]!}
|\Lambda;n_{lm}>=\prod_{l\leq m}e_q^{x_{lm}X_{lm}^-}|\Lambda>, \tag 31
$$
where $e_q^x$ is a $q$-exponential function:
$
e^x_q=\sum_{n=0}^\infty{x^n \over [n]!}
$.
The bilinear form $V_\Lambda^*
\otimes V_\Lambda \to {\Bbb C}$is defined by \quad$<\!\!\Lambda|\Lambda\!\!>=1$
and $(<\!\!{\bold x}|X)|{\bold y}\!\!>=<\!\!{\bold x})|(X|{\bold y}\!\!>)$
for any $<\!\!{\bold x}|\in V^*_\Lambda,|{\bold y}\!\!>\in V_\Lambda$
and $X\in \Uqslmn$. We can construct the representation in the space
of analytic functions on the coordinates $\{x_{ij}\}$ through
the bilinear form.
\par
\proclaim{Proposition 2}
The quantum superalgebra $\Uqslmn$ is realized by the Heisenberg algebra
${\Cal H}_L$.
$$
\align
t_i=&q^{-\sum^{M+N+1}_{m=1}\sum_{l=1}^m\sum_{r=l}^{m}a_{ir}M_{lm}
+(\alpha_i|\Lambda)}\tag 32\\
e_i=&q^{\sum_{k=1}^{i-1}(\nu_{i+1}M_{ki}-\nu_{i}M_{li-1})}{x_{ii}^{-1}}
[M_{ii}]+\sum_{l=1}^{i-1}
q^{\sum_{k=1}^{l-1}(\nu_{i+1}M_{ki}-\nu_{i}M_{li-1})}
x_{li-1}x^{-1}_{li}[M_{li}],\tag 33 \\
f_i=&\nu_{i}\sum_{j=1}^{i-1}q^{\rho_{ij}}
x_{ji}x_{ji-1}^{-1}[M_{ji-1}]
-\nu_{i+1}\sum_{j=i+1}^{M+N+1}q^{-\eta_{ij}}x_{ij}x_{i+1j}^{-1}[M_{i+1j}]  \\
&+x_{ii}\biggl[-{a_{ii} \over 2}M_{ii}-  
\sum^{M+N+1}_{m=i+1}\sum_{l=1}^m\sum_{r=l}^{m}a_{ir}M_{lm}
+(\alpha_i|\Lambda)\biggr],\tag 34 
\endalign
$$
where
$$
\align
\rho_{ij}&=-\sum^{M+N+1}_{m=i+1}\sum_{l=1}^m\sum_{r=l}^{m}a_{ir}M_{lm}
-\sum_{l=j+1}^i\sum_{r=l}^i a_{ir}M_{lj}
+\sum_{l=j+1}^{i-1}\nu_{i}M_{li-1}
+(\alpha_i|\Lambda), \tag 35\\
\eta_{ij}&=-\sum^{M+N+1}_{m=j+1}\sum_{l=1}^m\sum_{r=l}^{m}a_{ir}M_{lm}
-\sum_{l=i}^j\sum_{r=l}^j a_{ir}M_{lj}
+\sum_{l=i}^{j}a_{il}-\sum_{l=i+1}^{j}a_{il}
+(\alpha_i|\Lambda). \tag 36\\
\endalign
$$
\endproclaim 
\leftline{proof.}
\par
\noindent
We derive the forms by acting the generators $t_i,e_i,f_i(i=1,\cdots,M+N+1)$
on the $q$-deformed coherent states(31) and 
using the following commutation relations:
\footnote{Last three relations are derived from the Serre relations, the first two
use only one Serre relation,while the third one uses two Serre relations.
The detail calculations are given in the paper[22].}
$$
\align
&[h_i,e^{x_{lm}X_{lm}^-}_q]=-\sum_{r=l}^{m}a_{ir}M_{lm}e^{x_{lm}X_{lm}^-}_q,\tag 37 \\
&[e_i,e^{x_{jj}X_{jj}^-}_q]|\Lambda\!>=\delta_{i,j}
x_{ii}\biggl[-{a_{ii} \over 2}M_{ii}+(\alpha_i|\Lambda)\biggr]e^{x_{ii}X_{ii}^-}_q
|\Lambda\!>,\tag 38\\
&[e_i,X_{ji}^-]=\nu_{i}X_{ji-1}^-t_i, \tag 39\\
&[e_i,X_{ij}^-]=-\nu_{i+1}t_i^{-1}X_{i+1j}^-,\tag 40\\
&f_i(X_{ji-1}^-)^n=(X_{ji-1}^-q^{-\nu_i})^nf_i+[n](X_{ji-1})^{n-1} X_{ji}^-,\tag 41\\
&[f_i,X_{ji}^-]_{q^{\nu_i}}=0,\quad [f_j,X_{ji}^-]_{q^{-\nu_j}}=0,
\quad  [f_k,X_{ji}^-]=0  \quad (j+1\le k \le i-1). \tag 42\\
\endalign
$$
Substituting the explicit Cartan matrix $a_{ij}$,we obtain the following forms
in terms of $\nu_{i}$,$M_{ij}$ and $x_{ij}$:
$$
\align
\sum^{M+N+1}_{m=1}\sum_{l=1}^m\sum_{r=l}^{m}a_{ir}M_{lm}
&=\sum_{l=1}^{i-1}(\nu_{i+1}M_{li}-\nu_{i}M_{li-1}) \\
+\sum_{l=i+1}^{M+N+1}(\nu_{i}M_{il}-&\nu_{i+1}M_{i+1l})
+(\nu_i+\nu_{i+1})M_{ii}, \tag 43 \\
\sum^{M+N+1}_{m=i+1}\sum_{l=1}^m\sum_{r=l}^{m}a_{ir}M_{lm}
&=\sum_{l=i+1}^{M+N+1}(\nu_{i}M_{il}-\nu_{i+1}M_{i+1l}), \tag 44\\
\sum_{l=j+1}^{i}\sum_{r=l}^{i}a_{ir}M_{li}
&=\sum_{l=j+1}^{i-1}\nu_{i+1}M_{li}+(\nu_i+\nu_{i+1})M_{ii}, \tag 45\\
\sum_{l=i}^{j}\sum_{r=l}^{j}a_{ir}M_{lj}&=\nu_{i}M_{ij}-\nu_{i+1}M_{i+1j}, \tag 46\\
\sum_{l=i}^{j}a_{il}-\sum_{l=i+1}^{j}a_{il}&=\nu_i+\nu_{i+1}. \tag 47 \\
\endalign
$$
\proclaim{Proposition 3}
The superalgebra $\slmn$ is realized in terms of
$M_{ij}$, $x_{ij}$ ,$D_{ij}$ and $\nu_{ij}$.
$$
\align
t_i&=q^{-\sum_{l=1}^{i-1}(\nu_{i+1}M_{li}-\nu_{i}M_{li-1})
-\sum_{l=i+1}^{M+N+1}(\nu_{i}M_{il}-\nu_{i+1}M_{i+1l})
-(\nu_i+\nu_{i+1})M_{ii}+(\alpha_i|\Lambda) }, \tag 48 \\
e_{i}&=q^{\sum_{k=1}^{i-1}(\nu_{i+1}M_{ki}-\nu_{i}M_{ki-1})}D_{x_{ii}}
+\sum_{l=1}^{i-1}q^{\sum_{k=1}^{l-1}(\nu_{i+1}M_{ki}-\nu_{i}M_{ki-1})}
x_{li-1}D_{x_{li}}, \tag 49 \\
f_i&=\nu_{i}\sum_{j=1}^{i-1}q^{\rho_{ij}}
x_{ji}D_{x_{ji-1}}
-\nu_{i+1}\sum_{j=i+1}^{M+N+1}q^{-\eta_{ij}}x_{ij}D_{x_{i+1j}} \\
+&x_{ii}\biggl [-{\nu_i+\nu_{i+1} \over 2}M_{ii}-
\sum_{l=i+1}^{M+N+1}(\nu_{i}M_{il}-\nu_{i+1}M_{i+1l})
+(\alpha_i|\Lambda)\biggr],\tag 50 \\
\endalign
$$
where
$$
\align
\rho_{ij}&=-\sum_{l=j+1}^{i-1}(\nu_{i+1}M_{li}-\nu_{i}M_{li-1})
-\sum_{l=i+1}^{M+N+1}(\nu_{i}M_{il}-\nu_{i+1}M_{i+1l}) \\
&-(\nu_i+\nu_{i+1})M_{ii} 
+(\alpha_i|\Lambda), \tag 51\\
\eta_{ij}&=-\sum_{l=j}^{M+N+1}(\nu_{i}M_{il}-\nu_{i+1}M_{i+1l})
+(\nu_i+\nu_{i+1})+(\alpha_i|\Lambda). \tag 52\\
\endalign
$$
\endproclaim
\leftline{proof.}
\par
\noindent
It can be  checked that the generators satisfy the relations(20-26) 
by the explicit calculation of the commutation relations 
using the following relations:
$$
\align
[M_{ij}]x_{ij}&=x_{ij}[M_{ij}+1], \tag 53\\
[M_{ij}]D_{ij}&=D_{ij}[M_{ij}-1], \tag 54\\
[D_{ij},z_{lm}]&=\delta_{i,l}\delta_{j,m}\{[M_{lm}+1]-[M_{lm}]\}, \tag 55\\
[D_{ij},\theta_{lm}]&=\delta_{i,l}\delta_{j,m}\{[1-M_{lm}]+[M_{lm}]\}
=\delta_{i,l}\delta_{j,m}. \tag 56\\
\endalign
$$
We show the example of $U_q(sl(2|1)$.
$$
\align
t_1=&q^{-2M_{11}-M_{12}+M_{22}+(\alpha_1|\Lambda)}, \tag 57\\
t_2=&q^{M_{11}+M_{12}+(\alpha_2|\Lambda)}, \tag 58\\
e_1=&D_{x_{11}}, \tag 59\\
e_2=&q^{-M_{11}-M_{22}}\partial_{\theta_{22}}+x_{11}\partial_{\theta_{12}},  \tag 60\\
f_1=&-q^{M_{12}-M_{22}-(\alpha_1|\Lambda)-2}
\theta_{12}\partial_{\theta_{22}}+x_{11}[-M_{11}-M_{12}+M_{22}+(\alpha_1|\Lambda)], \tag 61\\
f_2=&q^{(\alpha_2|\Lambda)}\theta_{12}D_{x_{11}}+\theta_{22}[(\alpha_2|\Lambda)]. \tag 62\\
\endalign
$$
We give one concrete check.
$$
\align
[e_2,f_2]&=[q^{-M_{11}-M_{22}}\partial_{\theta_{22}},\theta_{22}[(\alpha_2|\Lambda)]]
+[x_{11}\partial_{\theta_{12}},q^{(\alpha_2|\Lambda)}\theta_{12}D_{x_{11}}] \\
&=q^{-M_{11}-M_{22}}[(\alpha_2|\Lambda)]+q^{(\alpha_2|\Lambda)}[M_{11}+M_{12}]
={t_2-t_2^{-1} \over q-q^{-1}}. \tag 63\\
\endalign
$$
\vfill\break
\enddocument